\documentclass{amsart}

\usepackage[latin1]{inputenc}
\usepackage{amssymb}
\usepackage{latexsym}
\usepackage{color}
\usepackage{longtable}

\newtheorem{teor}{Theorem}

\newtheorem{lema}[teor]{Lemma}
\newtheorem{prop}[teor]{Proposition}

\newcommand{\Q}{\mathbb{Q}}
\newcommand{\Qbar}{\overline{\mathbb{Q}}}
\newcommand{\Z}{\mathbb{Z}}
\newcommand{\F}{\mathbb{F}}

\newcommand{\Gal}{\operatorname{Gal\,}}

\newcommand{\GL}{\operatorname{GL}}

\newcommand{\Frob}{\operatorname{Frob}}

\newcommand{{\tors}}{\operatorname{tors}}

\def\QED{\hfill $\Box$}

\begin{document}

\title[On fields of definition of torsion points of elliptic curves with CM]{On fields of definition of torsion points of elliptic curves with complex multiplication}

\author[L. Dieulefait]{Luis Dieulefait}
\address{Universitat de Barcelona, Barcelona, Spain}
\email{ldieulefait@ub.edu}

\author[E.Gonz\'alez-Jim\'enez]{Enrique Gonz\'alez-Jim\'enez}
\address{Universidad Aut{\'o}noma de Madrid, Departamento de Matem{\'a}ticas and Instituto de Ciencias Matem{\'a}ticas (CSIC-UAM-UC3M-UCM), Madrid, Spain}
\email{enrique.gonzalez.jimenez@uam.es}

\author[J. Jim\'enez Urroz]{Jorge Jim\'enez Urroz}
\address{Universitat Polit{\'e}cnica de Catalunya, Barcelona, Spain}
\email{jjimenez@ma4.upc.edu }


\thanks{The authors were partially supported by the grants MTM2006-04895, MTM2009-07291 and CCG08--UAM/ESP--3906,  and MTM2009-11068 respectively.}

\date{\today}


\begin{abstract}
For any elliptic curve $E$ defined over the rationals with complex multiplication and for every prime $p$, we describe the image of the mod $p$ Galois representation attached to $E$. We deduce information about the field of definition of torsion points of these curves, in particular we classify all cases where there are torsion points over Galois number fields not containing the field of definition of the CM.
\end{abstract}
\maketitle

\section{Introduction}
Elliptic curves defined over the rationals have been extensively studied
for over a hundred years, mainly because we know that the set of
rational points form a finitely generated abelian group. In this paper we are interested in giving a
description of the possible finite groups that can be realized in this way. If we
focus on the torsion subgroup, the problem is solved:  we know every
possible group that appear as the rational torsion subgroup of an elliptic
curve over the rationals.

However, when considering more general number fields, many questions about
the torsion remain unsolved.  Here we approach this problem from the
following point of view:  we are interested in studying how does  the
torsion subgroup of an elliptic curve change when we enlarge the field of
definition. As a first step of this project, we are going to consider
elliptic curves defined over $\Q$.

The first results on this problem dealt with the case of quadratic number fields. For this case, Kwon \cite{kwon} has given results allowing to compute the torsion subgroup over a quadratic field of an elliptic curve defined over the rationals that has all the $2$-torsion subgroup defined over the rationals. First Qiu and Zhang \cite{QZ} and then completed by Fujita \cite{fujita} have generalized Kwon's result to the case of polyquadratic number fields.

We will focus on the complex multiplication (CM) case. We have proved results that give all the information about the field of definition of torsion points for CM elliptic curve defined over $\Q$. In particular, we have classified all cases where there are torsion points over Galois number fields not containing the field of definition of the CM. We also give a description of the image of the Galois representations on the $p$-torsion of all these curves for every prime $p$.

As it is well-known, the mod $p$ Galois representations attached to an elliptic curve behave in two very different ways depending on whether or not the curve has CM. Contrary to what happens in the CM case, in the non-CM case Serre \cite{serre72} proved that the images of the mod $p$ Galois representations are maximal, i.e., equal to $\GL_2(\F_p)$, for all but finitely many primes.

\

{\it \noindent Notation:} \\
$-$ Let $F$ be a number field, we will denote by $\mathcal{O}_F$ the ring of integers of $F$.\\
$-$ Let $p$ be an odd prime, we will denote by $\zeta_p$ a primitive $p^{\mbox{\tiny{th}}}$-root of the unity and by $\Q^+(\zeta_p)$ the maximal real subfield of the $p^{\mbox{\tiny th}}$-cyclotomic field $\Q(\zeta_p)$, that is $\Q^+(\zeta_p)=\Q(\zeta_p+\overline{\zeta}_p)$. Note that $\sqrt{-p}$ (respectively  $\sqrt{p}$) belongs to $\Q(\zeta_p)$ if $p\equiv\,3 (\mbox{mod}\,4)$ (respectively $p\equiv\,1 (\mbox{mod}\,4)$).\\
$-$  Let $E$ be an elliptic curve defined over $F$, the $m$-torsion subgroup of $E(F)$ will be denoted by $E(F)[m]=\{ P\in E(F)\,|\, [m]P=\mathcal{O} \}$ and by $E[m]=E(\overline{F})[m]$. \\
$-$ We will denote by $F(E[m])$ the number field obtained adjoining the coordinates of the points of order $m$.\\
$-$ Let $p$ be a prime, we will denote by $\rho_{E,p}$ the mod $p$ Galois representation attached to the $p$-torsion points of $E$.\\
$-$ $\chi$ will denote the mod $p$ cyclotomic character, where the
prime $p$ will be clear by the context.\\
$-$ We will denote by $j(E)$ the $j$-invariant of the elliptic curve $E$.

\

{\bf Acknowledgement:} We are grateful to John Cremona for useful comments at the beginning of this work. We also thank William Stein for computational facilities on the \verb+sage+ server. Finally, we thank the anonymous referee for useful comments.

\section{Elliptic curve with complex multiplication over $\Q$}

It is well known that there are 13 isomorphic classes of elliptic curves defined over $\Q$ with CM (cf. \cite[A \S 3]{advanced}). The following table gives a representative elliptic curve over $\Q$ for each class, that is, an elliptic curve over $\Q$ with CM by an order $R=\Z+\mathfrak{f}\,\mathcal{O}_K$ of conductor $\mathfrak{f}$ in a quadratic imaginary field $K=\Q(\sqrt{-D})$, where $ \mathcal{O}_K$ is the ring of integer of $K$. We will denote by $E_{D,\mathfrak{f}}$ this elliptic curve. 

\

\begin{center}
\begin{tabular}{|c|c|rcl|}
\hline
$-D$ & $\mathfrak{f}$ &&&{\!\!\!\!\!\!\!\!\!\!\!\!\!\!\! Short Weierstrass model of $E_{D,\mathfrak{f}}$}  \\
\hline
$-3$  & $1$ & $\quad y^2\!\!\!\!$ & $=$ & $\!\!\!\!x^3 + 16$\\
\hline
$-3$ & $2$ & $\quad y^2\!\!\!\!$ & $=$ & $\!\!\!\!x^3-15 x+22$\\
\hline
$-3$ & $3$ & $\quad y^2\!\!\!\!$ & $=$ & $\!\!\!\!x^3 - 480 x + 4048$\\
\hline
$-4$ & $1$ & $\quad y^2\!\!\!\!$ & $=$ & $\!\!\!\!x^3+x$\\
\hline
$-4$ & $2$ & $\quad y^2\!\!\!\!$ & $=$ & $\!\!\!\!x^3 - 11 x + 14$\\
\hline
$-7$ & $1$ & $\quad y^2\!\!\!\!$ & $=$ & $\!\!\!\!x^3 - 2835x - 71442$\\
\hline
$-7$ & $2$ & $\quad y^2\!\!\!\!$ & $=$ & $\!\!\!\!x^3-595x+5586$\\
\hline
$-8$ & $1$ & $\quad y^2\!\!\!\!$ & $=$ & $\!\!\!\!x^3 - 4320 x + 96768$\\
\hline
$-11$ & $1$ & $\quad y^2\!\!\!\!$ & $=$ & $\!\!\!\!x^3 - 9504 x + 365904$\\
\hline
$-19$ & $1$ & $\quad y^2\!\!\!\!$ & $=$ & $\!\!\!\!x^3 - 608 x + 5776 $\\
\hline
$-43$ & $1$ & $\quad y^2\!\!\!\!$ & $=$ & $\!\!\!\!x^3 - 13760 x + 621264$\\
\hline
$-67$ & $1$ & $\quad y^2\!\!\!\!$ & $=$ & $\!\!\!\!x^3 - 117920 x + 15585808 $\\
\hline
$-163$ & $1$ & $\quad y^2\!\!\!\!$ & $=$ & $\!\!\!\!x^3 - 34790720 x + 78984748304\quad$\\
\hline
\end{tabular}\\[5mm]
Table $1$ : Isomorphic classes of elliptic curves defined over $\Q$ with CM.
\end{center}

\

Let $E$ be an elliptic curves defined over $\Q$. We know that any $E'/\Q$ isomorphic to $E$ over $\overline{\Q}$ is in fact  $\Q$-isomorphic to a twist of $E$ (cf. \cite[X \S 5]{silverman}). More precisely, let $E:y^2=x^3+ax+b$ be a Weierstrass model for $E$, and $E'$ a curve isomorphic to $E$. We will
denote once and for all 
$$
n(E)=\left\{
\begin{array}{ccl}
2 & & \mbox{if $j(E)\ne 0,1728$,}\\
4 & & \mbox{if $j(E)= 1728$,}\\
6 & & \mbox{if $j(E)= 0$.}\\
\end{array}
\right.
$$
Then, $E'=E^d$ has  a Weierstrass model of the form
$$
\begin{array}{lclcl}
\mbox{(i)} & & E^d:y^2=x^3+d^2ax+d^3b & & \mbox{if $j(E)\ne 0,1728$,} \\
\mbox{(ii)} & & E^d:y^2=x^3+d ax & & \mbox{if $j(E) =1728$,}\\
\mbox{(iii)} & & E^d:y^2=x^3+d b & & \mbox{if $j(E) = 0$,}
\end{array}
$$ 
where $d$ is an integer in $\Q^*/(\Q^*)^{n(E)}$. In particular, any CM elliptic curve $E$ defined over  $\Q$  is in fact $\Q$-isomorphic to a curve $E^d_{D,\mathfrak{f}}$ for some $D,f$ as in Table $1$, and $d$ an integer in $\Q^*/(\Q^*)^{n(E)}$.

\

As we mentioned, in order to study the torsion of the elliptic curve, we will be using 
the mod $p$ Galois representation of the elliptic curve, for any prime $p$. It is then important to note
that, for $j(E)\ne 0, 1728$,  if $\rho_{E_{D,\mathfrak f},p}$ is the representation associated to $E_{D,\mathfrak f}$, then 
\begin{equation}\label{psi}
 \rho_{E_{D,\mathfrak f}^d,p}=\rho_{E_{D,\mathfrak f},p}\otimes\psi(d)
\end{equation}
is the twisted representation by the   Legendre symbol $\psi(d)=\left(\frac{d}{p}\right)$.  For the general case, we can  look at the number of points of the reduced curve to get information about the trace of the Frobenius at $p$, a spliting prime in the CM field. It is well known that the following formula holds  in general,  
$$
|E_{D,\mathfrak f}^d(\mathbb F_p)|=p+1+\pi\psi_{n(E)}(d)+\overline{\pi\psi_{n(E)}(d)},
$$
where $\pi$ is a primary prime above $p$ and $\psi_{n(E)}(\cdot)$ is the $n(E)$-power residue symbol.

\section{Statements of the main results}

\begin{teor}{($2$-torsion)}\label{teor2}
Let $E$ be an elliptic curve defined over $\Q$ with CM by an order of $K=\Q(\sqrt{-D})$ of conductor $\mathfrak{f}$ and let $F$ be a Galois number field not containing $K$, then
\begin{itemize}
\item $j(E)\ne 0,1728$:

\begin{itemize}
\item If $D\ne 8$ and $\mathfrak{f}$ odd, then $E(F)[2]=E(\Q)[2]$.
\item Otherwise, $\Q(E[2])=\Q(\sqrt{p})$ where $p|D$, in particular there are $2$-torsion points in a quadratic field different from $K$.
\end{itemize}
\item  $j(E)=1728$: In this case, $E=E_{4,1}^d$ for $d\in\Q^*/(\Q^*)^4$ and $\Q(E[2])=\Q(\sqrt{-d})$, in particular for $d\ne 1$ there are $2$-torsion points in a quadratic field different from $K$.
\item  $j(E)=0$: In this case, $E=E_{3,1}^d$ for $d\in\Q^*/(\Q^*)^6$ and $\Q(E[2])=\Q(\sqrt{-3},\sqrt[3]{2d})$. Moreover, $E(F)[2]=E(\Q)[2]$.
\end{itemize}
\end{teor}

\begin{teor}\label{teor:good}
Let $E$ be an elliptic curve defined over $\Q$ with CM by an order of $K=\Q(\sqrt{-D})$ and $p$ an  odd prime not dividing $D$. Let $F$ be a Galois number field not containing $K$, then $E(F)[p]$ is trivial.
\end{teor}

\begin{teor} \label{teor:bad}
Let $E$ be an elliptic curve defined over $\Q$ with CM by an order of $K=\Q(\sqrt{-D})$ of conductor $\mathfrak{f}$. We know that $E=E^d_{D,\mathfrak{f}}$ for some integer $d\in \Q^*/(\Q^*)^{n(E)}$. Let $p$ be an odd prime dividing $D$
\begin{itemize}
\item If $p>7$ then there are $p$-torsion points of $E$ defined over $\Q(\zeta_p+\overline{\zeta}_p,\sqrt{d})$. Furthermore,  $d=-p$ is the only case where any Galois number field containing $p$-torsion points contains $K$.

\item If $D=7$:
\begin{itemize}
\item Case $\mathfrak{f}=1$. There are $7$-torsion points of $E$ defined over $\Q(\zeta_7+\overline{\zeta}_7,\sqrt{-7d})$. Furthermore,  $d=1$ is the only case where any Galois number field containing $7$-torsion points contains $K$.
\item Case $\mathfrak{f}=2$. There are $7$-torsion points of $E$ defined over  $\Q(\zeta_7+\overline{\zeta}_7,\sqrt{7d})$. Furthermore,  $d=-1$ is the only case where any Galois number field containing $7$-torsion points contains $K$.
\end{itemize}

\item If $D=3$:
\begin{itemize}
\item  Case $\mathfrak{f}=1$. $\mathbb Q(E[3])=\mathbb Q(d^{1/6},\sqrt{-3})$. There is a  $3$-torsion point in the field $\Q(\sqrt{d})$ and, except for $d=-3$, this quadratic field is different from $K$. Moreover, 
if $d=e^3$,  there is a  $3$-torsion point on $\mathbb Q(\sqrt{-3e})$ which, except when $e$ is a square, is different from $K$.

\item Case $\mathfrak{f}\ne 1$. There are $3$-torsion points in the field $\Q(\sqrt{d})$. Except for $d=-3$ this quadratic field is different from $K$.
\end{itemize}
\end{itemize}
\end{teor}

{Remark: $p=3$ is the only odd prime where there are $p$-torsion points defined over $\Q$. For example in the last item of the previous theorem   $\mathfrak{f}=1$ with $e=1$ or $e=-3$, when $d=e^3$, and $\mathfrak{f}\ne 1$ with $d=1$ are the only cases where the curve has rational $3$-torsion.}

\section{On the division polynomials for small or bad reduction primes}

In this section we are going to study the $p$-division polynomial of the elliptic curves $E_{D,\mathfrak{f}}$, where $p=2,3$ or a bad reduction prime (i.e. $p=D$ with $2,3 \nmid D$). This study has been done computing explicitly the factorization of the $p$-division polynomial for the $13$ curves on the Table 1.

\

Let $E$ be an elliptic curve defined over a number field $F$ given by a short Weiestrass equation of the form $y^2=x^3+ax+b$, where $a,b\in\mathcal{O}_F$. Define the $m$-division polynomial $\Psi_m$, attached to $E$, recursively as follows:
$$
\begin{array}{l}
\psi_1 =1, \\
\psi_ 2= 2y\\
\psi_ 3= 3x^4+6ax^2+12bx-a^2,\\
\psi_ 4= (2x^6+10ax^4+40bx^3-10a^2x^2-8bax-2a^3-16b^2)\Psi_2\\
\psi_ {2k+1}=\psi_{k+2}\psi_{k}^3-\psi_{k-1}\psi_{k+1}^3,\,\,k\ge 2 \\
\psi_ {2k}=\psi_k (\psi_{k+2}\psi_{k-1}^2-\psi_{k-2}\psi_{k+1}^2)/\psi_{2},\,\,k\ge 2.
\end{array}
$$
Now, for $m>2$ define
$$
\Psi_m=\left\{
\begin{array}{lll}
\psi_{m} & & \mbox{$m$ odd}\\
\psi_{m}/\psi_2 & & \mbox{$m$ even}.
\end{array}
\right.
$$
and $\Psi_2=x^3+ax+b$. We have $\Psi_m\in\mathcal{O}_F[x]$. Then $P\in E[m]$ if and only if $\Psi_m(x(P))=0$. In particular,  $P\in E(F)[m]$ if and only if $\Psi_m(x(P))=0$ and $P\in E(F)$.

\

{\it Remark:} Let $E$ be an elliptic curve defined over a number field $F$ with $j(E)\ne 0,1728$ and $E^d$ its $d$-twist, for some $d\in F$. Then it is a straightforward computation to check that $\Psi_m^d(dx)=d^{\mbox{deg} \Psi_m} \Psi_m(x)$, where $\Psi_m$ (respectively $\Psi_m^d$) denotes the $m$-division polynomial of $E$ (respectively $E^d$). Therefore the study of the behavior of the roots of $\Psi^d_m(x)$ could be done on $\Psi_m(x)$ instead.

\

Let $E$ be an elliptic curve defined over $\Q$, $p$ a prime and $g(x)\in\Z[x]$ be a factor of the $p$-division polynomial of $E$. We denote by 
$$
\Q[E,g]=\Q(\{\alpha,\beta \in\Qbar\,|\, (\alpha,\beta)\in E[p]\,,g(\alpha)=0\}).
$$
Thus, $\Q(E[m])=\Q[E,\Psi_m]$.

\begin{lema}\label{lemma4}
Let $-D\in\{-7,-11,-19,-43,-67,-163\}$ and $E=E_{D,\mathfrak{f}}$ be an elliptic curve in the Table $1$. Let $p=D$ and $\Psi_p(x)$ be the $p$-division polynomial of $E$. Then the irreducible factorization of $\Psi_p(x)$ over $\Z[x]$ is given by
$$
\Psi_p(x)=g_p(x)h_p(x),\,\,\mbox{where}\,
\left\{
\begin{array}{l}
\mbox{deg}\,g_p=\,\,\,\,(p-1)/2\,,\\
\mbox{deg}\,h_p=p(p-1)/2\,.\\[2mm]
\end{array}
\right.
$$
Let us  denote by $F'_p=\Q(\{\alpha\in\Qbar \,|\,g_p(\alpha)=0\})$, then $F'_p=\Q^+(\zeta_p)$. Moreover, if we denote by $F_p=\Q[E,g_p]$ then :
\begin{enumerate}
\item[(i)] If $p\neq 7$, then $F_p=F_p'$. Therefore, $\Q(\zeta_p)=F'_p(\sqrt{-p})=F'_p\cdot \Q(\sqrt{-D})$.
\item[(ii)] If $p=7$ and $\mathfrak{f}=1$, then $F_7=\Q(\zeta_7)$. Thus, $\Q(\sqrt{-7})\subset F_7$.
\item[(iii)] If $p=7$ and $\mathfrak{f}=2$, then $F_7=\Q^+(\zeta_7)(\sqrt{7})$. Thus, $\Q(\zeta_7)\neq F_7$ and $\Q(\sqrt{-7})\not\subset F_7$.
\end{enumerate}
\end{lema}


\begin{lema}\label{lem:3-division}
Let $E=E_{D,\mathfrak{f}}$ be an elliptic curve in the Table $1$, $K=\Q(\sqrt{-D})$ and let $g(x)\in\Z[x]$ be a non-linear irreducible factor of the $3$-division polynomial of $E$. Then $K(\sqrt{-3})\subset \Q[E,g]$. Furthermore:
\begin{itemize}
\item[(i)] If $D= 3$ then $\Q(E_{D,\mathfrak{f}}[3])=K(\sqrt[3]{\mathfrak{f}})$.
\item[(ii)] If $D\ne 3$ then $\#{\rm Gal}(\Q(E_{D,\mathfrak{f}}[3])/\Q)=8$ or $16$.
\end{itemize}
%
%
\end{lema}

\begin{lema}\label{lema2torsion}
Let $E=E_{D,\mathfrak{f}}$ be an elliptic curve in the Table 1 and $K=\Q(\sqrt{-D})$, then:
\begin{enumerate}
\item[(i)] If $D\ne 8$ and $\mathfrak{f}$ odd, then $K\subset \Q(E[2])$.
\item[(ii)] Otherwise, $\Q(E[2])=\Q(\sqrt{p})$ where $p|D$.
\end{enumerate}

Moreover, the following table shows $\Q(E[2])$:

\

\begin{center}
\begin{tabular}{|c|c|l|}
\hline
$-D$ & $\mathfrak{f}$ & $\Q(E_{D,\mathfrak{f}}[2])$  \\
\hline
$-3$  & $1$ & $K(\sqrt[3]{2}) $\\
\hline
$-3$ & $2$ & $\Q(\sqrt{3}) $\\
\hline
$-3$ & $3$ & $K(\sqrt[3]{2}) $\\
\hline
$-4$ & $1$ & $K $\\
\hline
$-4$ & $2$ & $\Q(\sqrt{2}) $\\
\hline
$-7$ & $1$ & $K $\\
\hline
$-7$ & $2$ & $\Q(\sqrt{7}) $\\
\hline
$-8$ & $1$ & $\Q(\sqrt{2}) $\\
\hline
$-11$ & $1$ & $K(\alpha)$, $\alpha^3+\alpha^2+\alpha-1=0$\\
\hline
$-19$ & $1$ & $K(\alpha)$, $\alpha^3-\alpha^2+3\alpha-1=0$\\
\hline
$-43$ & $1$ & $K(\alpha)$,  $\alpha^3-3\alpha^2+7\alpha-1=0$\\
\hline
$-67$ & $1$ & $K(\alpha)$,  $\alpha^3-\alpha^2+7\alpha-9=0$\\
\hline
$-163$ & $1$ & $K(\alpha)$, $\alpha^3-9\alpha^2+85\alpha-227=0$\\
\hline
\end{tabular}
\end{center}
\end{lema}

\

\

The proof of the above lemmas is a straightforward computation, that has been done by using \verb+Magma+ and \verb+Sage+ (cf. respectively \cite{magma}, \cite{sage}). All the sources are available from \verb+http://www.uam.es/enrique.gonzalez.jimenez/research/tables/CM/+ .\\ Note that the case with more computer cost was the factorization of the $163$-division polynomial of $E_{163,1}$. This polynomial is of degree $13284$ with huge coefficients. The file that stores it has size around $280$ \verb+MB+. Then the factorization was done  using the functionality \verb+PARI+ on \verb+Sage+.

\section{$2$-Torsion}

{\it Proof of Theorem \ref{teor2}:} Let $E$ be an elliptic curve defined over $\Q$ with CM by an order $R=\Z+\mathfrak{f}\,\mathcal{O}_K$  in a imaginary quadratic field $K=\Q(\sqrt{-D})$. Then if $j(E)\ne 0,1728$, $E$ is $\Q$-isomorphic to $E^d_{D,\mathfrak{f}}$ for some square free integer $d$. We have $\Psi_2^d(dx)=d^3 \Psi_2(x)$, where $\Psi_2$ (respectively $\Psi_2^d$) denotes the $2$-division polynomial of $E_{D,\mathfrak{f}}$ (respectively $E^d_{D,\mathfrak{f}}$). That is, if $E_{D,\mathfrak{f}}\,:\,y^2=x^3+ax+b$ then  $\Psi_2(x)=x^3+ax+b$ and $\Psi_2^d(x)=x^3+d^2ax+d^3b$. Now, since the points of order $2$ are the ones that have ordinate zero we have that the field of definition of a point of order $2$ on $E^d_{D,\mathfrak{f}}$ is the same as the one on $E_{D,\mathfrak{f}}$. Now, thanks to the Lemma \ref{lema2torsion}, we have that if $F$ is a Galois field not containing $K$ and $D\ne 8$ and $\mathfrak{f}$ odd, then $E(F)[2]$ does not increase with respect to the $2$-torsion define $d$ over $\Q$. For the case $D=8$ or $\mathfrak{f}$ even we have that $\Q(E[2])=\Q(\sqrt{p})$ where $p|D$.

Now let $E$ be such that $j(E)=1728$, then $E=E_{4,1}^d:y^2=x^3 + dx$ for some fourth powerfree integer $d$. Then it is trivial to realize that $\Q(E[2])=\Q(\sqrt{-d})$. The case $j(E)=0$ is similar to the above case.

\QED

\section{Primes not dividing the discriminant}

\begin{teor} Let $E$ be an elliptic curve defined over $\Q$ with CM  by an order of $K=\Q(\sqrt{-D})$. Then, for every odd prime $p \nmid D$, the Galois representation corresponding to the $p$-torsion points of $E$ is irreducible.  
\end{teor}
 
{\it Proof:}  Since $E$ has CM, half of the traces of $\rho_{E,p}$ will be $0$, more precisely, for every prime $q$ inert on $K$ and of good reduction, the trace $a_q = 0$.

Suppose that $\rho_{E,p}$ is reducible. This means that, after semi-simplifying, we have  $\rho_{E,p} = \mu_1 \oplus \mu_2,$ with $\mu_1$, $\mu_2$ characters. It is well-known that these characters will take values on $\F_p^*$ (this follows from the fact that the Galois representations attached to an elliptic curve are odd). The determinant of $\rho_{E,p}$ is $\chi$, thus we can write $\mu_2 = \chi \cdot \mu_1^{-1}$.
 
Thus  $\rho_{E,p}$ is isomorphic to the sum of the two characters $\chi \mu_1^{-1}$ and $\mu_1$. Comparing traces for the image of  $\Frob \; q$ for any prime $q$ inert in $K$ and of good reduction for $E$ we obtain:
\begin{equation}\label{eqq}
a_q = 0  \equiv  q \mu_1(q)^{-1} + \mu_1(q) \qquad    \pmod{p}\,.
\end{equation} 
Since, by assumption, $p$ is not ramified in $K$, if we fix a non-zero residue class $t$ modulo $p$ Cebotarev density theorem implies that there are infinitely many primes $q$ that are inert in $K$ and also congruent to $t$ modulo $p$. We can also assume that these primes $q$ are of good reduction for $E$. We fix the following residue class modulo $p$: take $w$ to be any quadratic non-residue modulo $p$ and let $t = -w$. Thus, there are infinitely many primes $q$ inert in $K$, of good reduction for $E$, and congruent to $t= -w$ modulo $p$. For these primes $q$ congruence in (\ref{eqq}) gives:
$$ w \equiv -q \equiv \mu_1(q)^{2} \qquad    \pmod{p}\,.  $$
Since $-q \equiv w \; \pmod{p}$ which is a non-square this is a contradiction (recall that we know the character $\mu_1$ to take values on $\F_p^*$). This proves the theorem. \QED

\

{\it Proof of Theorem \ref{teor:good}:} We have shown that the representation is irreducible. On the other hand, by the theory of complex multiplication the restriction to  ${\rm Gal}(\overline{\Q}/K)$ has abelian image. Therefore the image  contains an abelian normal subgroup of index $2$, it has a dihedral  proyectivization like
\begin{equation}
\left<\left( \begin{array}{cc}
 * & \, \\
\, & * \end{array}\right),\left( \begin{array}{cc}
 0 & 1 \\
1 & 0 \end{array}\right)  \right>,
\end{equation}
where the element $T$ of order $2$ comes from the conjugation 
$\sigma\in\Gal(K/\Q)$. Hence, for any field $K\subsetneq F\subset E[p]$, $\Gal(F/\Q)$ will be a quotient of 
$\Gal(E[p]/\Q)$ containing $T$ and, in particular, again irreducible. The conclusion follows since it is obvious that a torsion point defined over $F$ will produce a subspace invariant by  $\Gal(F/\Q)$ which would then give a reducible representation.\QED

\

The following proposition describes the image of the mod $p$ Galois representation obtained in this section.

\begin{prop}
Let $E$ be an elliptic curve defined over $\Q$ with CM by an order of $K=\Q(\sqrt{-D})$ and $p$ an  odd prime not dividing $D$. Then the proyectivization of $\rho_{E,p}$ has dihedral image.
\end{prop}

\section{Primes dividing the discriminant}
{\it Proof of Theorem \ref{teor:bad}:}
We consider a curve $E=E_{D,\mathfrak{f}}$ as in Table $1$ and $p=D > 3$ which
is a bad reduction prime. We will use Lemma \ref{lemma4} together
with the information that the curve has CM defined over $K$ to give
a precise description of the Galois number field generated by the
$p$-torsion points of $E$. \\
The image of $\rho_{E,p}$ is
contained in $\GL_2(\F_p)$ and since the curve has CM over $K$, it
is well-known that the restriction to the Galois group of $K$ is
reducible since the $p$-torsion points generate an abelian extension
of $K$. Thus, the image of $\rho_{E,p}$ is a group containing an abelian normal subgroup with index at most $2$.\\
 We know by Lemma \ref{lemma4} that the order
of the Galois number field generated by the $p$-torsion of $E$ is
divisible by $p$ (in fact, a factor $h_p$ of the $p$-division
polynomial has degree divisible by $p$). This Galois number field is
the one corresponding to the image of $\rho_{E,p}$. If we
apply Dickson's classification (cf. \cite[Prop. 15]{serre72}) of maximal subgroups of $\GL_2(\F_p)$
we conclude that the representation is either reducible or
surjective (we know that the determinant is surjective). But since
we also know that it contains an abelian normal subgroup of index at
most $2$ it has to be reducible. Since its order is divisible by $p$
it is reducible but not decomposable. Thus, the image of $\rho_{E,p}$
can be described as follows: 
\begin{equation}\label{repr}
\rho_{E,p} \cong  \left(
\begin{array}{ccc}
          \phi_1 &   *\\
          0 &   \phi_2 \end{array} \right)
\end{equation}
          where $* \neq 0$ and $\phi_1, \phi_2$ are characters of
          $\Gal(\overline{\Q}/\Q)$ satisfying $\phi_1 \cdot \phi_2 = \chi$. It is known that these two characters must have
           values in $\F_p^*$, because the representation $\rho_{E,p}$ is odd and therefore if it is reducible
            it must reduce over $\F_p$. In this
          case, precisely because the representation is reducible,
          either there are torsion points defined over $\Q$ (but it
          is well-known that there are no such points since $p>3$) or
          there is a non-trivial Galois invariant $1$-dimensional subspace, the
          character $\phi_1$ corresponding to the action of $\Gal(\overline{\Q}/\Q)$
          on this subspace. In Lemma \ref{lemma4} we observe that
          there is an abelian extension generated by some of the
          torsion points (the one corresponding to points of torsion whose $x$-coordinates are roots of $g_p$), this abelian extension must clearly
          correspond to the Galois invariant subspace with character
          $\phi_1$. In Lemma \ref{lemma4} this is the extension that we have called $F_p$.  Then we can recover the character $\phi_1$ which corresponds to $F_p$,
          and also $\phi_2 = \chi \cdot \phi_1^{-1}$, in each case. From Lemma \ref{lemma4}, the character $\phi_1$ is: 
          \begin{itemize}
          \item[] if $ D=p > 7$, prime: $\phi_1=\chi^2$.
          \item[] if $ D=7$, $\mathfrak{f}=1$: $\phi_1=\chi$.
          \item[] if $ D=7$, $\mathfrak{f}=2$: $\phi_1=\mu$, where $\mu$ corresponds to the extension $\Q(\zeta_7+\overline{\zeta}_7,\sqrt{7})$.
          \end{itemize}
          
Let us now consider the case $p=3$, $\mathfrak{f}>1$ and $E=E_{3,\mathfrak{f}}$. These elliptic curves have rational $3$-torsion points, and therefore the representation $\rho_{E,3}$ is reducible as in (\ref{repr}) with $\phi_1=1$ and $\phi_2=\chi$ but,  by Lemma \ref{lem:3-division}, $\rho_{E,3}$ do not decompose as the sum of these two characters. For the general case of a CM curve $E$, we will use again the fact that $E=E^d_{D,\mathfrak{f}}$ and $\rho_{E,p}=\rho_{E_{D,\mathfrak{f}},p}\otimes \psi$ where $\psi$ is quadratic. Therefore we have for the image of $\rho_{E,p}$ that again it is contained in a Borel as in (\ref{repr}):
\begin{equation}\label{repr2}
\rho_{E,p} \cong  \left(
\begin{array}{ccc}
          \phi_1\cdot \psi &   *\\
          0 &   \phi_2\cdot \psi \end{array} \right)
\end{equation}
where  $*\ne 0$.

\

For $p=3$, $\mathfrak f=1$, we just have to consider the $3$-division polynomial of $E^d_{3,\mathfrak{1}}$, $\psi_3(x)=3x(x^3+2^6d)$,
to conclude that the coordinates of the $3$-torsion points are $(0,\pm 4\sqrt d)$, $(4\alpha d^{1/3},\pm 4\sqrt{-3d})$, where $\alpha$ is any cubic root of $-1$. From here,  (\ref{repr}) and (\ref{repr2}), Theorem \ref{teor:bad} follows easily . \QED

\

The description of the image of the mod $p$ Galois representation just obtained can be summarized as follows:

\begin{prop}
Let $E=E^d_{D,\mathfrak{f}}$ and $p=D$ be a prime greater than $3$. Then the image of $\rho_{E,p}$ is:
$$
\rho_{E,p} \cong  \left(
\begin{array}{ccc}
          \phi_1\cdot \psi &   *\\
          0 &   \phi_2\cdot \psi \end{array} \right)
$$
where 
$$
\phi_1=\left\{
\begin{array}{ccl}
\chi & & \mbox{if $ p=7$ and $\mathfrak{f}=1$}\\
\mu & & \mbox{if $ p=7$ and $\mathfrak{f}=2$}\\
\chi^2 &&  \mbox{if $p>7$}
\end{array}
\right.
$$
where $\mu$ corresponds to the extension $\Q(\zeta_7+\overline{\zeta}_7,\sqrt{7})$, $\phi_2=\chi\cdot\phi_1^{-1}$, $\psi$ is the quadratic character of $\mathbb Q(\sqrt d)$ and $* \neq 0$.

\end{prop}



\end{document}